\documentstyle[multicol,11pt]{article}
\begin{document}

\centerline{\Large\bf{On Robust Stability of Multivariable Interval
Control Systems}
}
\vspace*{2\baselineskip}

\centerline{\large Zhizhen Wang$^{(1,2)}$ \ \ \ Long Wang$^{(3)}$ \\
Wensheng Yu$^{(2)}$ \ \ \ Guoping Liu$^{(2,4)}$}
\vspace*{2\baselineskip}
\centerline{(1)\ \ Shanghai Normal University, Shanghai, 200234, CHINA}
\centerline{(2)\ \ Laboratory of Complex
Systems and Intelligence Science, Institute of Automation,}
\centerline{ Chinese Academy of Sciences, Beijing, 100080, CHINA}
\centerline{(3)\ \ Center for Systems and Control, Department of
Mechanics and Engineering Science,}
 \centerline{Peking University,
Beijing, 100871, CHINA }
\centerline{(4)\ \ School of Mechanical,
Materials, Manufacturing Engineering and Management, }
\centerline{University of Nottingham, United Kingdom }

\abovedisplayskip=0.13cm
\abovedisplayshortskip=0.06cm
\belowdisplayskip=0.13cm
\belowdisplayshortskip=0.08cm
\renewcommand{\baselinestretch}{1}

\vspace*{4\baselineskip}

{\small  {\bf Abstract:} This paper studies robustness of MIMO
control systems with parametric uncertainties, and establishes a
lower dimensional robust stability criterion. For  control systems
with interval transfer matrices, we identify the minimal testing
set whose stability can guarantee the stability of the entire
uncertain set. Our results improve the results in the literature,
and provide a constructive solution to the robustness of a family
of MIMO control systems.}

\vspace*{1\baselineskip} {\small {\bf Keywords:} Robust Stability,
Multivariable Feedback Systems, Polynomial Matrices, Interval
Polynomials, Edge Theorem.}

\section{Introduction}
\par\indent

Robustness is considered as an elementary property in the analysis
and design of control systems. A control system is said to be
robust if it retains an assigned degree of stability or
performance under perturbations. Robust control under structured
perturbations has been an active research area in recent
years$^{\cite{khar}-\cite{barm1}}$.

Robust stability problem for control systems with interval
characteristic polynomials was initiated by
Kharitonov(1978)$^{\cite{khar}}$, who showed that the Hurwitz
stability of a real interval polynomial family can be guaranteed
by the Hurwitz stability of four prescribed critical vertex
polynomials in this family. Bartlett, Hollot and Huang extended
the result and presented the Edge Theorem$^{\cite{bart}}$, which
says that the stability of a polytope of polynomials can be
guaranteed by the stability of its one-dimensional exposed edge
polynomials. Chapellat et al.(1989)$^{\cite{chap1}}$ refined the
Kharitonov's Theorem, and proved that the robust stability of the
unity feedback control system family with an interval plant and a
fixed controller in forward path can be guaranteed by the
stability of 32 prescribed one-dimensional unity feedback control
systems in this family.

When considering multivariable control systems, the robust
stability problem becomes more complicated and challenging
$^{\cite{bias}-\cite{wang}}$. For Hurwitz stability of interval
matrices, Bialas 'proved' that in order to guarantee robust
stability, it suffices to check all vertex
matrices$^{\cite{bias}}$. Later, it was shown by Barmish that
Bialas's result was incorrect${\cite{barm}}$. Kokame and Mori
established a Kharitonov-like result on robust Hurwitz stability
of interval polynomial matrices$^{\cite{koka}}$. Their result is
refined by Wang, Wang and Yu$^{\cite{long}}$. Kamal and Dahleh
established some robust stability criteria for MIMO control
systems with fixed controllers and uncertain
plant$^{\cite{kama}}$.

In this paper, we focus on a class of MIMO feedback control
systems with their transfer function matrices described by
\begin{equation}\label{1}
{\cal M}(s)={\cal B}(s)A(s)+{\cal D}(s)C(s)
\end{equation}
\begin{equation}\label{2}
{\cal B}_{ij}(s)=\left\{
\begin{array}{l}
b_{ij}^{(0)}+b_{ij}^{(1)}s+\cdots+b_{ij}^{(n)} s^n, \  \
b_{ij}^{(k)}\in [\underline b_{ij}^{(k)},\overline b_{ij}^{(k)}]
\end{array}\right\}
\end{equation}
\begin{equation}\label{3}
{\cal D}_{ij}(s)=\left\{
\begin{array}{l}
d_{ij}^{(0)}+d_{ij}^{(1)}s+\cdots+d_{ij}^{(n)}
s^n, \  \
d_{ij}^{(k)}\in[\underline d_{ij}^{(k)},\overline d_{ij}^{(k)}]
\end{array}\right\}
\end{equation}
\begin{equation}\label{4}
{\cal B}(s)=\left\{\left(b_{ij}(s)\right)_{n\times n},\; b_{ij}(s)\in {\cal B}_{ij}(s)\right\}
\end{equation}
\begin{equation}\label{5}
{\cal D}(s)=\left\{\left(d_{ij}(s)\right)_{n\times n},\; d_{ij}(s)\in {\cal D}_{ij}(s)\right\}
\end{equation}
where $n$ is a given positive integer, and $\underline
b_{ij}^{(k)},\overline b_{ij}^{(k)}, \underline
d_{ij}^{(k)},\overline d_{ij}^{(k)}$ are given real numbers. And
$A(s)=\left(a_{ij}(s)\right)_{n\times
n},C(s)=\left(c_{ij}(s)\right)_{n\times n}$ are two fixed polynomial
matrices, ${\cal B}(s),{\cal D}(s)$ are interval polynomial
matrices and ${\cal B}_{ij}(s),{\cal D}_{ij}(s)$ are interval
polynomials, where $A(s)=\left(a_{ij}(s)\right)_{n\times n}$ means
that the matrix is an $n\times n$ one and its $ij$-th entry is
$a_{ij}(s)$. Rewrite ${\cal M}(s)$ as
$$
{\cal M}(s)=M_n s^n+M_{n-1} s^{n-1}+\cdots+M_0
$$
where $M_i$ is an $n\times n$ scalar matrix.

{\bf Assumption A:} $\mbox{det}(M_n)\not=0$ over ${\cal M}(s)$.

In the sequel, we study the multivariable interval control system
whose characteristic polynomial matrix ${\cal M}(s)$ satisfies
Assumption A. By exploiting the uncertainty structures, we are
able to reduce the computational burden in checking robust
stability.

An outline of the paper is as follows. Some useful definitions and
lemmas are presented in Section 2. Section 3 contains our main
result on minimal dimension robust stability criterion. An example
is worked out in Section 4, and some discussions are made in
Section 5. Concluding remarks are given in Section 6.

\section{Preliminaries}
\par\indent

{\bf Definition 2.1}\ \ Given an interval polynomial set ${\cal
R}(s)=\sum_{i=0}^m r_i s^i, \ \  r_i \in [r_i^L, r_i^U] $, its
Kharitonov vertex set is defined as
 ${\cal R}^V(s)=\{r_1(s), r_2(s), r_3(s), r_4(s)\}$, and its Kharitonov edge set
 is defined as
 ${\cal R}^E(s)=\{\lambda r_i(s)+(1-\lambda) r_j(s), \ \
(i,j)\in \{(1,2), (2,4), (4,3), (3,1)\}, \ \  \lambda \in [0,1]\}$,
 where
  $$
  \begin{array}{l}
  r_1(s)=r_0^L+r_1^L s+r_2^U s^2+r_3^U s^3+r_4^L s^4+\cdots\\
  r_2(s)=r_0^L+r_1^U s+r_2^U s^2+r_3^L s^3+r_4^L s^4+\cdots\\
  r_3(s)=r_0^U+r_1^L s+r_2^L s^2+r_3^U s^3+r_4^U s^4+\cdots\\
  r_4(s)=r_0^U+r_1^U s+r_2^L s^2+r_3^L s^3+r_4^U s^4+\cdots
  \end{array}
  $$
Denote  the Kharitonov vertex sets and the Kharitonov edge sets of
${\cal B}_{ij}(s),{\cal D}_{ij}(s)$ as ${\cal B}_{ij}^V(s),{\cal
D}_{ij}^V(s),{\cal B}_{ij}^E(s)$ and ${\cal D}_{ij}^E(s)$,
respectively.

{\bf Definition 2.2}\ \ A fix polynomial is said to be (Hurwitz)
stable if all its roots lie in the strict left half of the complex
plane. A polynomial matrix is said to be stable if its determinant
is stable.

{\bf Lemma 1}\quad (Edge Theorem$^{\cite{bart}}$) Let $\Omega$ be a polytope of polynomials without degree dropping.
Then the boundary of $R(\Omega)$ is contained in the root space of the exposed edges of $\Omega$,
where $R(W)=\{s \mbox{ such that } T(s)=0 \mbox{ for some polynomial } T\in W\}$ is called the root space of $W$.

{\bf Lemma 2}\ \  $A(s)$ and $C(s)$  are defined as before. Let $B(s)=\left(b_{ij}(s)\right)_{n\times n},
D(s)=\left(d_{ij}(s)\right)_{n\times n}$ and
\begin{equation}
\begin{array}{l}
N^{(i)}(s)=
\left\{
B(s)A(s)+D(s)C(s):
\begin{array}{l}
\exists k_1,k_2\in\{1,\cdots,n\};\ \
 l=1,\cdots,n\\
b_{ik_1}(s)\in {\cal B}^E_{ik_1}(s),\ \
 b_{il}(s)\in {\cal B}^V_{il}(s), \ \l\not=k_1;\\
d_{ik_2}(s)\in {\cal D}^E_{ik_2}(s),\ \
 d_{il}(s)\in {\cal D}^V_{il}(s),\ \ l\not=k_2
\end{array}\right\}
\end{array}
\end{equation}
\begin{equation}
\begin{array}{l}
N^{(i)}_E(s)=
\left\{B(s)A(s)+D(s)C(s):
\begin{array}{l}
\exists k\in\{1,\cdots,n\};\ \ l\not=k,\ \  l=1,\cdots,n\\
b_{ik}(s)\times d_{ik}(s)
\in \begin{array}{l}
\left({\cal B}^E_{ik}(s)\times {\cal D}^V_{ik}(s)\right)\cup \left({\cal B}^V_{ik}(s)\times {\cal D}^E_{ik}(s)\right)
\end{array}\\
b_{il}(s)\times d_{il}(s)
\in  \left({\cal B}^V_{il}(s)\times {\cal D}^V_{il}(s)\right)
\end{array}\right\}
\end{array}
\end{equation}
Then, $N^{(i)}(s)$  is stable for all $i\in\{1,\cdots,n\}$ if and only if  $N^{(i)}_E(s)$ is stable.\\
Proof\ \ (Sufficiency) Suppose that $N^{(i)}_E(s)\mbox{ is stable }$, our aim is to prove that
$N^{(i)}(s) \mbox{ is stable}$.
For any $N(s)\in N^{(i)}(s)$, there exist $k_1,k_2\in\{1,\cdots,n\}$ such that
{\small
\begin{equation}\label{10}
 \begin{array}{ll}
 b_{ik_1}(s)\in {\cal B}^E_{ik_1}(s);& b_{il}(s)\in {\cal B}^V_{il}(s),\; l\not=k_1,l=1,\cdots,n\\
d_{ik_2}(s)\in {\cal D}^E_{ik_2}(s);& d_{il}(s)\in {\cal D}^V_{il}(s),\;
l\not=k_2,l=1,\cdots,n
\end{array}
\end{equation}
}
{
\arraycolsep=1mm
\begin{equation}\label{11}
\begin{array}{l}
N(s)=
\left(\begin{array}{ccc}
\sum_k \left(
\begin{array}{l}
b_{1k}(s)a_{k1}(s)
+d_{1k}(s)c_{k1}(s)
\end{array}\right)&\cdots&\sum_k
\left(\begin{array}{l}
b_{1k}(s)a_{kn}(s)
+d_{1k}(s)c_{kn}(s)
\end{array}\right)\\
\cdots&\cdots&\cdots\\
\sum_k \left(\begin{array}{l}
b_{ik}(s)a_{k1}(s)
+d_{ik}(s)c_{k1}(s)
\end{array}\right)&\cdots&\sum_k
\left(\begin{array}{l}
b_{ik}(s)a_{kn}(s)
+d_{ik}(s)c_{kn}(s)
\end{array}\right)\\
\cdots&\cdots&\cdots\\
\sum_k \left(
\begin{array}{l}
b_{nk}(s)a_{k1}(s)
+d_{nk}(s)c_{k1}(s)
\end{array}\right)&\cdots&\sum_k
\left(\begin{array}{l}
b_{nk}(s)a_{kn}(s)
+d_{nk}(s)c_{kn}(s)
\end{array}\right)
\end{array}
\right)
\end{array}
\end{equation}
}
By definition, $N(s)$ is stable if and only if ${\mathrm{det}}N(s)$ is stable. By Laplace formula, expanding the determinant
along the $i$-th row, we have
{
$$
\begin{array}{ll}
{\mathrm{det}}N(s)=&\sum_k \left(b_{ik}(s)a_{k1}(s)+d_{ik}(s)c_{k1}(s)\right)N_{i1}+\cdots
+\sum_k \left(b_{ik}(s)a_{kn}(s)+d_{ik}(s)c_{kn}(s)\right)N_{in},
\end{array}
$$
} where $N_{ij}$ is the algebraic complementary minor of the
$ij$-th entry of $N(s)$, which is independent of
$b_{i1}(s),\cdots, b_{in}(s),d_{i1}(s),\cdots,d_{in}(s)$. By a
simple algebraic manipulation, we have
$$
\begin{array}{rl}
{\mathrm{det}}N(s)=&\sum_{j=1}^n b_{ij}(s)\left(\sum_{k=1}^n
a_{jk}(s)N_{ik}\right)+\sum_{j=1}^n d_{ij}(s)\left(\sum_{k=1}^n c_{jk}(s)N_{ik}\right)\\
=& b_{ik_1}(s)\delta_1(s)+d_{ik_2}(s)\delta_2(s)+\delta_3(s)
\end{array}
$$
where $\delta_l(s)$ is a term which is independent
of $b_{ik_1}(s),d_{ik_2}(s)$. Denote
{
\begin{equation}\label{14}
\left\{\begin{array}{l} b_{ik_1}(s)\times d_{ik_2}(s)\in
\left(\begin{array}{l}
\left({\cal B}^E_{i k_1}(s)\times {\cal D}^V_{i k_2}(s)\right)
\cup \left({\cal B}^V_{i k_1}(s)\times {\cal D}^E_{i k_2}(s)\right)
\end{array}\right)\\
b_{il}(s)\in {\cal B}^V_{il}(s),\ \ l\not=k_1,\ \ l=1,\cdots,n\\
d_{il}(s)\in {\cal D}^V_{il}(s),\ \ l\not=k_2,\ \ l=1,\cdots,n
\end{array}\right.
\end{equation}
} By Lemma 1, the stability of each such polytope is equivalent to
the stability of its exposed edges. These edges are obtained by
letting $B(s)$ and $D(s)$ satisfy (\ref{14}). Thus, ${\mathrm{det}}N(s)$
is stable for $N(s)$ satisfying (\ref{10}),(\ref{11}) if and only
if ${\mathrm{det}}N(s)$ is stable for $N(s)$ satisfying
(\ref{11}),(\ref{14}). Since $N^{(i)}_E(s)\mbox{ is stable }$,
this shows that
$N(s)$ is stable. Hence, sufficiency is proved. \\
(Necessary) This is obvious since $N^{(i)}_E(s)\subset
N^{(i)}(s)$.

\section{ Robust Stability Criterion}
\par\indent

In this section, we consider the uncertain system family
(\ref{1}). Before presenting our main result, we first define some
notations. The notation $X\times Y$ stands for the Cartesian
product, which is the set of all ordered pairs $(x,y)$, where
$x\in X,y\in Y$. Let $S_n$ denote the set of all bijections of
the set $\{1,\cdots,n\}$ onto itself, and $T_n$ be the set of all functions from
the set $\{1,\cdots,n\}$ onto itself. Take
$$
\begin{array}{l}
{\cal B}_E(s)=
\left\{B(s)\in{\cal B}(s):
\begin{array}{l}
b_{ij}(s)\in {\cal B}^E_{ij}(s),  j=\sigma(i)\\
b_{ij}(s)\in {\cal B}^V_{ij}(s),  j\not=\sigma(i)\\
\end{array},\ \  \sigma\in S_n \right\}\\
{\cal D}_E(s)=
\left\{D(s)\in{\cal D}(s):
\begin{array}{l}
d_{ij}(s)\in {\cal D}_{ij}^E(s),  j=\sigma(i)\\
d_{ij}(s)\in {\cal D}_{ij}^V(s),  j\not=\sigma(i)\\
\end{array},\ \ \sigma\in S_n\right\}\\
{\cal B}_E{\cal D}_E={\cal B}_E(s)\times {\cal D}_E(s)\\
{\cal BD}=
\left\{ \begin{array}{l}
B(s)\times D(s)\\
\in {\cal B}_E{\cal D}_E
\end{array}:
\begin{array}{l}
\eta\in T_n\\
b_{ij}(s)\times d_{ij}(s)
\in \begin{array}{l}
\left({\cal B}^E_{ij}(s)\times {\cal
D}^V_{ij}(s)\right)
\cup \left({\cal B}^V_{ij}(s)\times
{\cal D}^E_{ij}(s)\right)\end{array},\;
 j=\eta(i);\\
b_{ij}(s)\times d_{ij}(s)
\in {\cal B}^V_{ij}(s)\times {\cal D}^V_{ij}(s),\;
j\not=\eta(i)\\
\end{array}\right\}
\end{array}
$$
The following theorem was given by Kamal and Dahleh (1996).

{\bf Proposition 1}$^{\cite{kama}}$\ \ ${\cal M}(s)$ is stable for
all $B(s)\in{\cal B}(s),D(s)\in{\cal D}(s)$ if and only if ${\cal
M}(s)$ is stable for all $B(s)\times D(s)\in {\cal B}_E{\cal
D}_E$.

As a generalization of  the Generalized Kharitonov Theorem,
Proposition 1 addresses the robust stability of ${\cal M}(s)$ and
reduces it to $4^{2n^2}(n!)^2$ robust stability problems involving
$2n$ uncertain parameters. In this paper, we will show that, by
further making use of the uncertainty structure information, the
original robust stability problem can be reduced to $(2n)!
4^{2n^2}(n!)$ robust stability problems involving $n$ uncertain
parameters. By Proposition 1 and Lemma 2, we get the main result
of this paper.

{\bf Theorem 1}\ \ ${\cal M}(s)$ is stable for all $B(s)\in{\cal
B}(s),D(s)\in{\cal D}(s)$ if and only if
${\cal M}(s)$ is stable for all $B(s)\times D(s)\in {\cal BD}$.

Proof:\ \ By Proposition 1, ${\cal M}(s)$ is stable for all
$B(s)\in{\cal B}(s),D(s)\in{\cal D}(s)$ if and only if ${\cal
M}(s)$ is stable for all $B(s)\times D(s)\in  {\cal B}_E{\cal
D}_E$. For any $n(s)\in {\cal M}(s)$ with $B(s)\times D(s)\in
{\cal B}_E{\cal D}_E$, fixing the rows from the second to the
$n$-th rows and applying Lemma 2 to the first row, we have
$$
\mbox{${\cal M}(s)$ is stable for all $B(s)\in{\cal
B}(s),D(s)\in{\cal D}(s)$ }
$$
if and only if
$$
\mbox{${\cal M}(s)$ is stable for all $B(s)\times D(s)\in {\cal BD}_1\cap {\cal B}_E{\cal D}_E$}
$$
where ${\cal BD}_1\subset {\cal B}(s)\times {\cal D}(s)$ satisfies
(\ref{14}) when $i=1$. Repeating the procedure for the remaining
rows consecutively, we get the conclusions.

{\bf Remark 1}\ \ Theorem 1 shows that, the criterion for the
robust stability of interval polynomial matrices can be greatly
simplified. Proposition 1 of Kamal  and Dahleh, reduced the
$2n^2$-dimensional stability verification problem to
$2n$-dimensional problem, and our result reduces the
$2n^2$-dimensional stability verification problem to
$n$-dimensional one.

{\bf Remark 2}\ \ Theorem 1  establishes a minimal dimension
testing criterion for robust stability of multivariable control
systems. By making use of the uncertainty information of  ${\cal
B}(s)$ and ${\cal D}(s)$, i.e., not only making use of the
uncertainty structure information in ${\cal B}(s)$ and ${\cal
D}(s)$ separately (this is exactly what Kamal and Dehleh did), but
also making use of the uncertainty structure information
collectively, we are able to improve the robust stability  results
given by Kamal and Dehleh.

{\bf Remark 3}\ \ Theorem 1 allows us to improve the main results
of Kamal and Dehleh.

\section{Illustrative Example}
\par\indent

Recall the two-link planar manipulator considered by Kamal and Dehleh.
The characteristic polynomial of the closed loop system is
\begin{equation}\label{ex}
T(\epsilon)=\det \left(M(\theta_d)s^3+K_ds^2+K_p s+K_r\right)
\end{equation}
where $\theta_d\in[0,\pi/2]$ is the joint angle and
$$
\begin{array}{l}
M(\theta_d)=\left(
\begin{array}{cc}
3+2\cos(\theta_d)&1+\cos(\theta_d)\\
1+\cos(\theta_d)&1
\end{array}\right)\\
K_d=\left(\begin{array}{cc}
k_{d11}&k_{d12}\\
k_{d21}&k_{d22}\end{array}\right)\\
K_p=\left(\begin{array}{cc}
k_{p11}&k_{p12}\\
k_{p21}&k_{p22}
\end{array}\right)\\
K_r=\left(\begin{array}{cc}
k_{r11}&k_{r12}\\
k_{r21}&k_{r22}
\end{array}\right)
\end{array}
$$
\begin{equation}\label{in}
\begin{array}{c}
k_{d11}\in[6.07-6.07\epsilon,6.07+6.07\epsilon],\ \
k_{d12}\in[2.22-2.22\epsilon,2.22+2.22\epsilon]\\
k_{d21}\in[2.22-2.22\epsilon,2.22+2.22\epsilon],\ \
k_{d22}\in[1.62-1.62\epsilon,1.62+1.62\epsilon]\\
k_{p11}\in[6.12-6.12\epsilon,6.12+6.12\epsilon],\ \
k_{p12}\in[2.24-2.24\epsilon,2.24+2.24\epsilon]\\
k_{p21}\in[2.24-2.24\epsilon,2.24+2.24\epsilon],\ \
k_{p22}\in[1.64-1.64\epsilon,1.64+1.64\epsilon]\\
k_{r11}\in[5.11-5.11\epsilon,5.11+5.11\epsilon],\ \
k_{r12}\in[1.87-1.87\epsilon,1.87+1.87\epsilon]\\
k_{r21}\in[1.87-1.87\epsilon,1.87+1.87\epsilon],\ \
k_{r22}\in[1.37-1.37\epsilon,1.37+1.37\epsilon]
\end{array}.
\end{equation}
Since for any fixed $\theta_d\in[0,\pi/2]$, the inertia matrix
{
$$
M(\theta_d)=
\left(
\arraycolsep=1mm
\begin{array}{cc}
1&(1+\cos(\theta_d))\\
(-1+\cos(\theta_d))&1
\end{array}\right)
\left(
\arraycolsep=1mm
\begin{array}{cc}
1&0\\
2&1
\end{array}\right)
$$}
In what follows, we take the uncertain inertia matrix as
$$
M(s)=B(s)A(s)
$$
where $A(s)=\left(\begin{array}{ll}1&0\\
2&1
\end{array}\right)$ and
$$B(s)\in{\cal B}(s)=\{(b_{ij}(s))_{2\times 2}:b_{ij}(s)\in {\cal B}_{ij}(s)\}$$
$$
\begin{array}{lll}
{\cal B}_{11}(s)=s^3 ,& {\cal B}_{12}(s)=b_{12}s^3,& b_{12}\in[1,2]\\
{\cal B}_{22}(s)=s^3,  & {\cal B}_{21}(s)=b_{21}s^3,& b_{21}\in[-1,0]
\end{array}
$$
Let ${\cal D}(s)=K_d s^3+K_ps^2+K_r$, then
$$
{\cal D}(s)=\left\{(d_{ij}(s)):d_{ij}(s)\in {\cal D}_{ij}(s)\right\},
$$
where
$$
{\cal D}_{ij}(s)=\left\{d_{ij}(s)=k_{dij}s^2+k_{pij}s+k_{rij}\right\}
$$
and $k_{dij},k_{pij} \mbox{ and } k_{rij}\mbox{ are given in (\ref{in})}$.
Thus, the characteristic polynomial of the uncertain system is
\begin{equation}
\begin{array}{l}
T(\epsilon)=\det\left(B(s)A(s)+D(s)\right); \ \
B(s)\in {\cal B}(s),\ \
D(s)\in {\cal D}(s)
\end{array}
\end{equation}
where $A(s)$ is a fixed matrix and ${\cal B}(s),{\cal D}(s)$ are interval polynomial matrices.

The Kharitonov polynomial vertex sets associated with ${\cal B}_{ij}(s)$
and ${\cal D}_{ij}(s)$ are as follows,
$$
\arraycolsep=0.6mm
\begin{array}{ll}
{\cal B}^V_{11}(s)=\{s^3\}={\cal B}_{11}(s), & {\cal B}^V_{12}(s)=\{s^3,2s^3\}, \\
 {\cal B}^V_{21}(s)=\{-s^3,0\}, & {\cal B}^V_{22}(s)=\{s^3\}={\cal B}_{22}(s),
\end{array}
$$
$$
\begin{array}{c}
{\cal D}^V_{ij}(s)=\left\{\begin{array}{c}
d_{ij}^{(1)}=K_{dij}^U s^2+K_{pij}^L s+K_{rij}^L,\ \
d_{ij}^{(2)}=K_{dij}^U s^2+K_{pij}^U s+K_{rij}^L,\\
d_{ij}^{(3)}=K_{dij}^L s^2+K_{pij}^L s+K_{rij}^U,\ \
d_{ij}^{(4)}=K_{dij}^L s^2+K_{pij}^U s+K_{rij}^U,
\end{array}\right\}
\end{array}
$$
where the upper bound and the lower bound of $K_{dij}$ are denoted as $K_{dij}^U,K_{dij}^L$.
And the Kharitonov edge sets associated with ${\cal B}(ij)$ and ${\cal D}(ij)$ are  as follows.
$$
\begin{array}{l}
{\cal B}^E_{11}(s)=\{s^3\}={\cal B}_{11}(s), \\
{\cal B}^E_{12}(s)=\{\lambda s^3+(1-\lambda)\cdot 2s^3:\lambda\in[0,1]\},\\
{\cal B}^E_{21}(s)=\{-\lambda s^3:\lambda\in[0,1]\}, \\
{\cal B}^E_{22}(s)=\{s^3\}={\cal B}_{22}(s)
\end{array}
$$
$$
{\cal D}^E_{ij}(s)=\left\{
\begin{array}{c}
\lambda d_{ij}^{(1)}+(1-\lambda)d_{ij}^{(2)},\ \
\lambda d_{ij}^{(2)}+(1-\lambda)d_{ij}^{(4)}\\
\lambda d_{ij}^{(4)}+(1-\lambda)d_{ij}^{(3)},\ \
\lambda d_{ij}^{(3)}+(1-\lambda)d_{ij}^{(1)}
\end{array}\right\}
$$
By Proposition 1 given by Kamal and Dehleh, for any $\epsilon$, $T(\epsilon)$ is robustly stable
if and only if $T(\epsilon)$ is stable for all $B(s)\times D(s)\in {\cal B}_E{\cal D}_E$, where
$$
{\cal B}_E {\cal D}_E=U_1\cup U_2
$$
{
$$
\begin{array}{c}
U_1= \left(\begin{array}{cc}
{\cal B}_{11}(s)&{\cal B}^E_{12}(s)\\
{\cal B}^E_{21}(s)&{\cal B}_{22}(s)
\end{array}\right)\times
\left(\begin{array}{cc}
{\cal D}^V_{11}(s)&{\cal D}^E_{12}(s)\\
{\cal D}^E_{21}(s)&{\cal D}^V_{22}(s)
\end{array}\right)\\
U_2= \left(\begin{array}{cc}
{\cal B}_{11}(s)&{\cal B}^E_{12}(s)\\
{\cal B}^E_{21}(s)&{\cal B}_{22}(s)
\end{array}\right)\times
\left(\begin{array}{cc}
{\cal D}^E_{11}(s)&{\cal D}^V_{12}(s)\\
{\cal D}^V_{21}(s)&{\cal D}^E_{22}(s)
\end{array}\right)
\end{array}
$$}
So, the robust stability of $T(\epsilon)$ is reduced to $2\cdot
4^6$ robust stability problems involving four parameters. In this
case, the testing set is a four-dimensional set. In the sequel, we
apply Theorem 1 to simplify this problem further. That is to say,
$T(\epsilon)$ is robustly stable if and only if $T(\epsilon)$ is
stable for all $B(s)\times D(s)\in {\cal B}{\cal D}$, where
{
$$
{\cal B}{\cal D}=V_1\cup V_2\cup V_3\cup V_4\cup V_5\cup V_6\cup V_7
$$
$$
\begin{array}{c}
V_1=\left(\begin{array}{cc}
{\cal B}_{11}(s)&{\cal B}^E_{12}(s)\\
{\cal B}^E_{21}(s)&{\cal B}_{22}(s)
\end{array}\right)\times
\left(\begin{array}{cc}
{\cal D}^V_{11}(s)&{\cal D}^V_{12}(s)\\
{\cal D}^V_{21}(s)&{\cal D}^V_{22}(s)
\end{array}\right)\\
V_2=\left(\begin{array}{cc}
{\cal B}_{11}(s)&{\cal B}^E_{12}(s)\\
{\cal B}^V_{21}(s)&{\cal B}_{22}(s)
\end{array}\right)\times
\left(\begin{array}{cc}
{\cal D}^V_{11}(s)&{\cal D}^V_{12}(s)\\
{\cal D}^V_{21}(s)&{\cal D}^E_{22}(s)
\end{array}\right)\\
V_3=\left(\begin{array}{cc}
{\cal B}_{11}(s)&{\cal B}^E_{12}(s)\\
{\cal B}^V_{21}(s)&{\cal B}_{22}(s)
\end{array}\right)\times
\left(\begin{array}{cc}
{\cal D}^V_{11}(s)&{\cal D}^V_{12}(s)\\
{\cal D}^E_{21}(s)&{\cal D}^V_{22}(s)
\end{array}\right)\\
V_4=\left(\begin{array}{cc}
{\cal B}_{11}(s)&{\cal B}^V_{12}(s)\\
{\cal B}^E_{21}(s)&{\cal B}_{22}(s)
\end{array}\right)\times
\left(\begin{array}{cc}
{\cal D}^V_{11}(s)&{\cal D}^E_{12}(s)\\
{\cal D}^V_{21}(s)&{\cal D}^V_{22}(s)
\end{array}\right)\\
V_5=\left(\begin{array}{cc}
{\cal B}_{11}(s)&{\cal B}^V_{12}(s)\\
{\cal B}^E_{21}(s)&{\cal B}_{22}(s)
\end{array}\right)\times
\left(\begin{array}{cc}
{\cal D}^E_{11}(s)&{\cal D}^V_{12}(s)\\
{\cal D}^V_{21}(s)&{\cal D}^V_{22}(s)
\end{array}\right)\\
V_6=\left(\begin{array}{cc}
{\cal B}_{11}(s)&{\cal B}^V_{12}(s)\\
{\cal B}^V_{21}(s)&{\cal B}_{22}(s)
\end{array}\right)\times
\left(\begin{array}{cc}
{\cal D}^E_{11}(s)&{\cal D}^V_{12}(s)\\
{\cal D}^V_{21}(s)&{\cal D}^E_{22}(s)
\end{array}\right)\\
V_7=\left(\begin{array}{cc}
{\cal B}_{11}(s)&{\cal B}^V_{12}(s)\\
{\cal B}^V_{21}(s)&{\cal B}_{22}(s)
\end{array}\right)\times
\left(\begin{array}{cc}
{\cal D}^V_{11}(s)&{\cal D}^E_{12}(s)\\
{\cal D}^E_{21}(s)&{\cal D}^V_{22}(s)
\end{array}\right)
\end{array}
$$
}
Clearly, the robust stability of $T(\epsilon)$ is reduced to $7\cdot 4^6$ robust stability problems involving
two parameters. In this case, the minimal testing set is two-dimensional.

{\bf Remark 4:}\ \ A test with one uncertain parameter could
involve infinite tests, e.g., an edge contains infinite points.
The significance of our result is that it is the minimal dimension
test for such problem (\ref{1}). Although the number of robust
stability problem is slightly increased, the dimension of the
testing set is significantly reduced. Furthermore, we considered
the uncertainty in inertia matrix $M(\theta_d)$ in the above
Example, whereas Kamal and Dahleh did not take into account
\cite{kama}.

\section{Discussions}
\par\indent

The focal point of this paper is to establish a lower dimensional robust stability criterion for robustness of
multivariable control systems. By improving the result given by Kamal and Dahleh, we have significantly reduced
the high-dimensional stability verification problem to a much lower-dimensional problem in parameter space,
e.g., for a $2\times 2$ interval control system family, our result reduces $8$-dimensional
robust problem to a $2$-dimensional one and Dehleh's result reduced this an $8$-dimensional
robust problem to a $4$-dimensional one. Furthermore, when dealing with SISO systems, our result naturally reduces to
the well-known Box Theorem. Of course, from the computational viewpoint, checking robust $D$-stability
of interval polynomial matrices is NP-hard. Hence, some sufficient (conservative) conditions were
developed based on linear matrix inequality(LMI) techniques. And how to reduce the conservativeness
of these conditions is still under investigations. on the other hand, from the theoretical viewpoint, establishing
lower-dimensional necessary and sufficient condition for robust $D$-stability of interval polynomial
matrices is important as well. By our result, the 'nucleus' of the uncertainty structures can be identified, which
is very useful in dealing with robust performance of control systems under structured parametric
perturbations$^{\cite{khar}-\cite{wang1}}$.

\section{Conclusions}
\par\indent

In this paper, robustness of MIMO control systems with interval
transfer function matrices is studied. For a family of
multivariable control systems having interval transfer functions
matrices, a lower dimensional robust stability criterion has been
established. This result improves the previous criterion obtained
by Kamal and Dahleh, and can be used to establish improved lower
dimensional robust performance criteria as well.

\section*{Acknowledgements}
\par\indent

This work was supported by National Natural Science Foundation of
China (10372002, 60204006 and 69925307), the Knowledge Innovation
Pilot Program of the Chinese Academy of Sciences, National Key
Basic Research and Development Program (No. 2002CB312200).

\end{document}